\documentclass[12pt]{article}
\usepackage{amsmath, amsthm, amssymb}

\begin{document}
\begin{center}
{\large \bf Repetition in Colored Sequences of Balls}\\
Jeremy M. Dover
\end{center}
\begin{abstract}
In responding to a question on Math Stackexchange, the author~\cite{2109013} formulated the problem of determining the number of sequences of balls colored in most $n$ colors where some fixed numbers of balls share colors with other balls in the sequence. In this paper, we formulate the problem more formally, and solve several variations.
\end{abstract}

\section{The setup}
Suppose you have balls of $n$ different colors, with an ample supply of each color, and you wish to create a sequence of $k$ balls. We assume that balls of the same color are indistinguishable, but positions in the sequence are distinguishable. A very natural question is to ask how many such sequences contain exactly $m$ balls which are the same color as another ball. Conceptually, this problem is easily solved in terms of sums of multinomial coefficents, but the challenge is figuring out which selections of colors to sum over. To break the problem down more effectively, we pose a more specific problem: how many ways are there to color a sequence of $k$ balls with at most $n$ colors such that exactly $m$ balls are the same color as another ball, and exactly $\lambda$ colors are repeated. We denote this number as $Z(k,n,m,\lambda)$.

First, note some constraints, added to the obvious constraints that all of $k$, $n$, $m$, and $\lambda$ must be non-negative integers:
\begin{itemize}
\item If $k>n$, then there must be at least $k-n+1$ balls that match some other ball, since at most $n$ balls are the first instance of their color in the sequence, meaning the remaining $k-n$ balls must match one of these, and at least one of these first instances must be matched. Hence $Z(k,n,m,\lambda) = 0$ for all $m<k-n+1$.
\item If exactly $m$ balls match some other ball, then at most $\lfloor \frac{m}{2} \rfloor$ colors can be repeated. Hence $Z(k,n,m,\lambda) = 0$ for $\lambda > \lfloor \frac{m}{2} \rfloor$.
\item If exactly $m$ balls match some other color, then $k-m$ balls do not match any other color. Hence the number of colors not repeated in the sequence, $n - \lambda$, must be at least $k-m$. Hence $Z(k,n,m,\lambda) = 0$ for all $\lambda > n-k+m$.
\item If $m=0$, then $\lambda$ is necessarily 0 as well, and we have $Z(k,n,0,0) = \frac{k!}{(k-n)!}$ for all $n \le k$, and $Z(k,n,0,\lambda) = 0$ for all other choices of $k$, $n$, and $\lambda$.
\item We cannot have exactly one ball matching, so $Z(k,n,1,\lambda) = 0$ for all $k$, $n$, and $\lambda$
\end{itemize}

Given that we have values for $k$, $n$, $m$, and $\lambda$ which meet the constraints above, the count is fairly straightforward. First, select the $\lambda$ colors to be repeated, which can be done in $n \choose \lambda$ ways. Next, select the $m$ balls in the sequence which are going to match other positions, which can be done in $k \choose m$ ways. Now color the $k-m$ balls which need to be uniquely colored with the remaining $n-\lambda$ colors, which can be done in $\frac{(n-\lambda)!}{(n-\lambda+m-k)!}$ ways.

The final component is to color the $m$ balls that will match with the $\lambda$ repeated colors such that each color is used at least twice. Thinking of the color assignment as a function, this function is doubly-surjective, that is, each color is the image of at least two balls in the domain.  Doubly-surjective functions have been studied by Walsh~\cite{walsh}, who defines $s(m,\lambda)$ to be the number of doubly-surjective functions from a set of $m$ elements into a set of $\lambda$ elements. Walsh calculates the following formula:
$$s(m,\lambda) = \sum_{j=0}^\lambda (-1)^j{\lambda \choose j}\sum_{i=0}^j {j \choose i}\frac{m!}{(m-i)!}(\lambda-j)^{m-i}$$

With this notation, we have $$Z(k,n,m,\lambda) = {n \choose \lambda}{k \choose m}\frac{(n-\lambda)!}{(n-\lambda+m-k)!}s(m,\lambda)$$ where $m\ge k-n+1$ and $\lambda \le {\rm min}\{\lfloor \frac{m}{2} \rfloor, n-k+m\}$.

\section{Counting matching balls}
With the formula in the previous section, we have the building blocks to solve several similar, but slightly different problems.

{\em Problem 1: How many sequences of $k$ balls colored in at most $n$ colors contain exactly $m \ge 2$ balls with a color matching some other ball?}

 From our previous analysis, any sequence satisfying the conditions of the problem must have $\lambda$ colors which are repeated, where $1 \le \lambda \le {\rm min}\{\lfloor \frac{m}{2} \rfloor, n-k+m\}$. Therefore the solution to this problem is simply

$$\sum_{\lambda=0}^{{\rm min}\{\lfloor \frac{m}{2} \rfloor, n-k+m\}} Z(k,n,m,\lambda)$$

As an example, we can count the sequences of 5 balls colored in at most three colors such that exactly 4 balls match some other color. In this case we could have $\lambda=1$ or $\lambda=2$. If $\lambda=1$ we must have four balls of one color, which can be chosen in three ways, and one ball of a differing color, which can be chosen in two ways. Once the balls are chosen, there are five ways to put them in sequence, namely $AAAAB, AAABA, \ldots, BAAAA$. There are 30 such sequences.

The case for $\lambda=2$ is slightly more complex. Since four balls match and two colors are repeated, we must have two balls of each of two colors, and one of the third color; these choices can be made in three ways, since all colors are effectively assigned once the singleton color is picked. Once the balls are picked, there are five ways to place the singleton, then ${4 \choose 2} = 6$ ways to place the first matching pair. This yields $3 \times 5 \times 6 = 90$ possible sequences of this form. Summing with the previous case gives 120 possible solutions, which matches the formula calculation.

{\em Problem 2: How many sequences of balls colored in at most $n$ colors contain exactly $m \ge 2$ balls with a color matching some other ball?}

This problem is similar to the previous, but note that the length constraint has been removed. For fixed $m$ and $n$, the only sequences colored with at most $n$ colors and having exactly $m$ balls matching another ball must have length between $m$ and $m+n-1$, using our constraints above. This can be answered via summing over the appropriate answers to the previous problem:

$$\sum_{k=m}^{m+n-1}\sum_{\lambda=0}^{{\rm min}\{\lfloor \frac{m}{2} \rfloor, n-k+m\}} Z(k,n,m,\lambda)$$

\section{Counting repeats}

Conceptually, when counting matches in the previous problem, we assume that the entire sequence has already been created, and we look at each ball of the sequence and can determine if some other ball has a matching color, no matter where in the sequence the match arises. Instead, suppose we are drawing the balls sequentially, and must make a ``matching'' decision when the ball is drawn; thus we are counting the number of balls in the sequence whose color has previously been seen in the sequence. This produces a subtle difference in our counting, because we do not go back and retrospectively count the first occurrence of a color that will later match.

To see the difference, consider the sequence $AABBCCDDDD$. This sequence has 10 balls and four colors, and all 10 balls match some other color. However, only 6 balls match a color previously seen in the sequence, since the first instance of each of the four colors is not counted.

{\em Problem 3: How many sequences of $k$ balls colored in at most $n$ colors contain exactly $\mu \ge 1$ balls whose color appears earlier in the sequence?}

As before, there are some constraints on values for $k$, $n$ and $\mu$ for this problem to make sense, beyond the obvious non-negativity constraints. First note that $\mu \le k-1$, since the first ball of a sequence will never be counted, though all of the remaining balls may be. In addition, only the first instance of a color in a sequence is not counted as a repeat, and there are at most $n$ such first instances, so we must have $k \le n+\mu$.

Using the same trick as before, let $\lambda$ be the number of colors which occur at least twice in a particular sequence. Clearly, we must have $1 \le \lambda \le \mu$, and any sequence with exactly $\mu$ balls whose color appears previously in the sequence has exactly $\mu+\lambda$ balls whose color matches some other in the sequence ($\mu$ repeats, plus the original occurrences of each of $\lambda$ colors). Therefore we can use the numbers $Z(k,n,m,\lambda)$ to solve this problem as well, with the answer being
$$\sum_{\lambda=0}^\mu Z(k,n,\mu+\lambda,\lambda)$$

{\em Problem 4: How many sequences of balls colored in at most $n$ colors contain exactly $\mu \ge 1$ balls whose color appears earlier in the sequence?}

Like the problems in the previous section, for fixed $\mu$ and $n$ there only exist sequences of balls colored with at most $n$ colors and having exactlu $\mu$ balls match a color previously seen for lengths $k$ satisfying $\mu +1 \le k \le n+\mu$. As before, we can simply sum the corresponding values for the individual lengths to obtain the answer:

$$\sum_{k=\mu+1}^{n+\mu}\sum_{\lambda=0}^\mu Z(k,n,\mu+\lambda,\lambda)$$

\bibliography{bibfile}{}
\bibliographystyle{plain}
\end{document}